\documentclass[11pt]{article}
\usepackage{amsfonts,amssymb,amsmath,amscd}
\usepackage{hyperref}
 \usepackage{authblk}
\usepackage{color}
\usepackage{graphicx}
\usepackage{wrapfig}   
\usepackage{changepage}
\newcommand{\R}{\mathbb{R}}

\newcommand{\bfx}{{\bf x}}
\newcommand{\bfy}{{\bf y}}

\newcommand{\beq}{\begin{equation}}
\newcommand{\eeq}{\end{equation}}
\newcommand{\beqs}{\begin{eqnarray}}
\newcommand{\eeqs}{\end{eqnarray}}

\title{
Robert V. Kohn (1953-2026)
}

\begin{document}

\author[1]{Yekaterina Epshteyn
  \thanks{\href{epshteyn@math.utah.edu}{epshteyn@math.utah.edu}}}
\author[2]{David Kinderlehrer
  \thanks{\href{davidk@andrew.cmu.edu}{davidk@andrew.cmu.edu}}}
\affil[1]{\small Department of Mathematics, The University of Utah, UT, USA}
\affil[2]{\small Department of Mathematical Sciences, Carnegie Mellon
  University, Pittsburgh, PA, USA}

\maketitle

\begin{abstract}
 The tribute article is dedicated to the memory and enduring legacy of Professor Robert
 V. Kohn, Courant Institute, NYU. In the article we record thoughts and reminiscences of his exemplary life.

\end{abstract}
\section*{Introduction}

Robert Vita Kohn, Bob Kohn to all, born October 5, 1953, in Shaker Heights, Ohio, passed away at home in New York, January 12, 2026, at the age of 72 from cancer.
\begin{figure}
\centering{\includegraphics[width=0.52\textwidth]{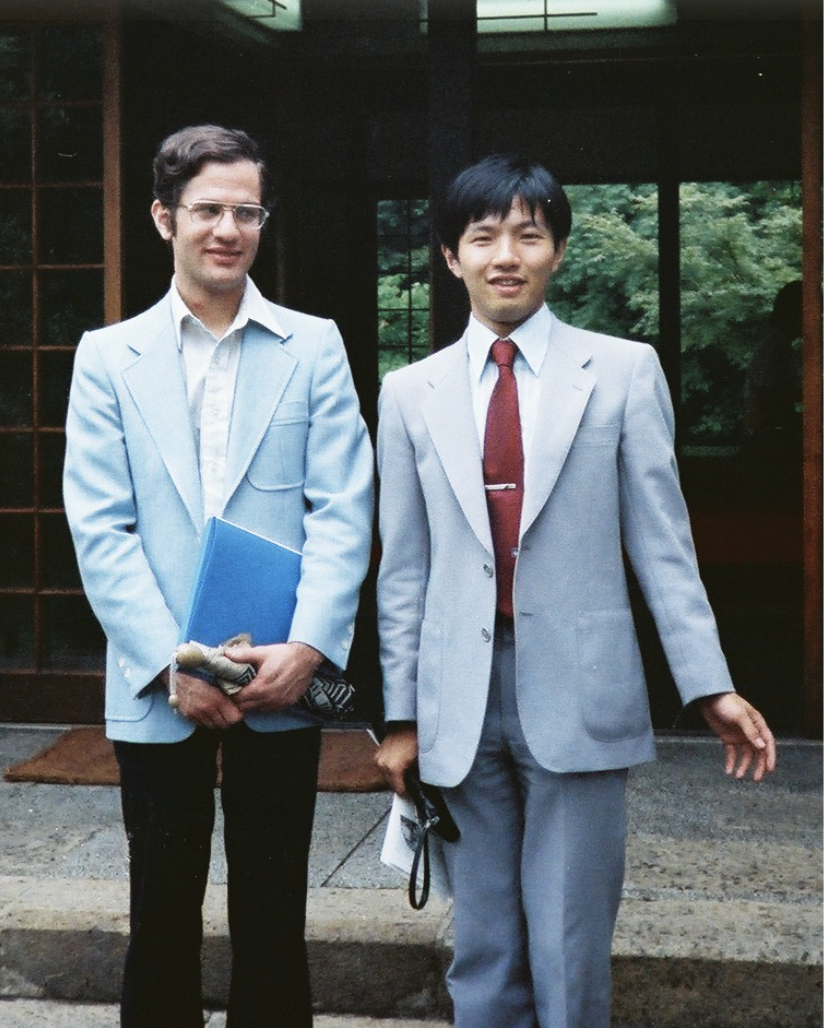}
    \quad
\includegraphics[width=0.52\textwidth]{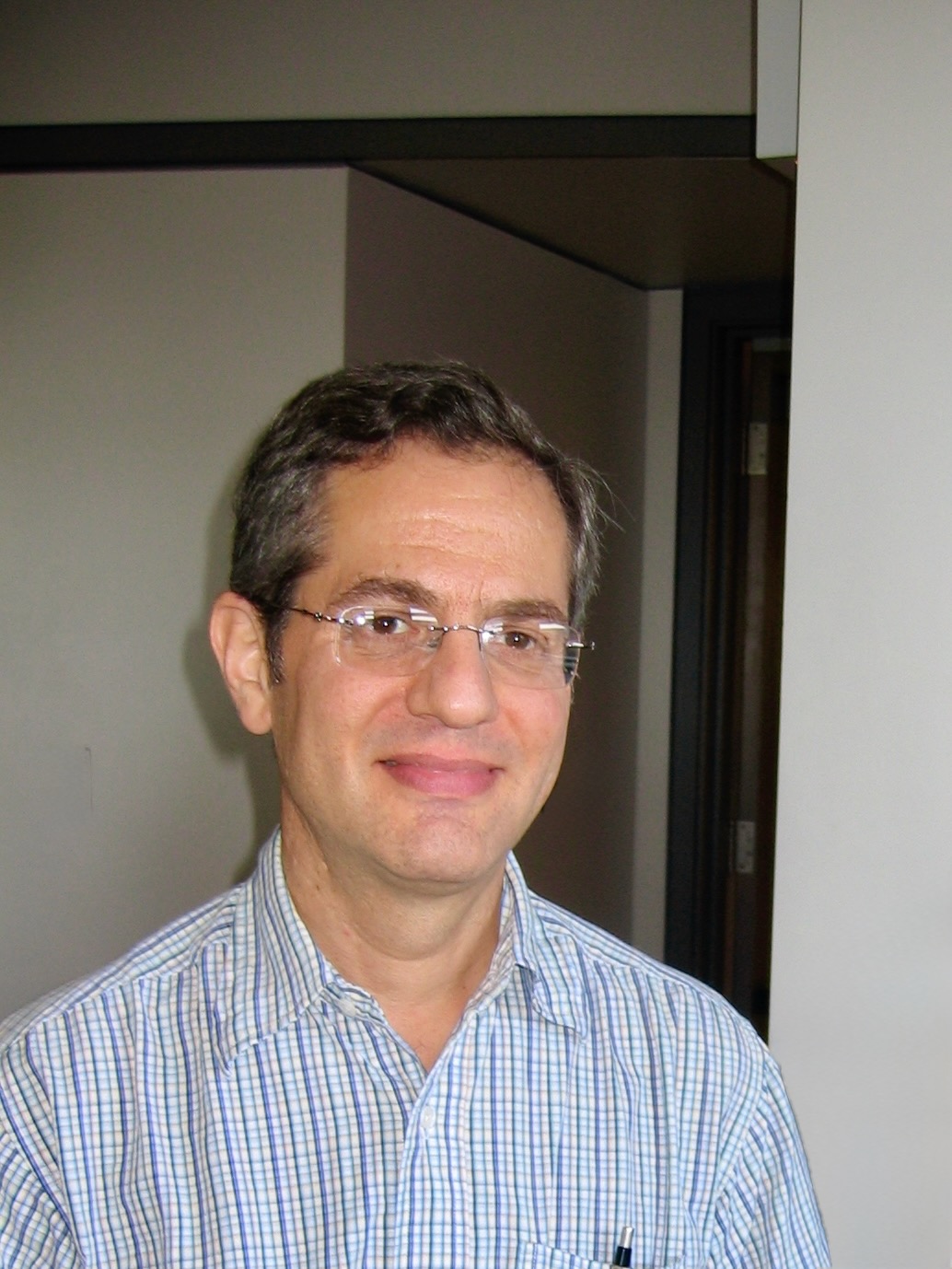}}
\caption{\footnotesize R.V. Kohn and Y. Giga, 1982,  Nagoya (top) - photo due to Y. Giga. R.V. Kohn, 2006, Montreal (bottom) - photo due to D. Kinderlehrer.
} \label{fig1}.
\end{figure}
The Courant Institute has a storied legacy, now enhanced by Bob's profound scientific vision, his mentoring, and his service.  In this Memorial article we record thoughts and reminiscences of this exemplary life.

Bob was raised in a loving family with his siblings, Luci and Norman.  He was fascinated with math and science from an early age.  When in high school, Bob attended and later served as a counselor during summers at the Ohio State University mathematics school.  After he completed high school a year early, he studied at Harvard University, earning an A.B. in Mathematics in 1974.
He took a Masters at Warwick, publishing his first papers in logic. Bob Kohn received a PhD in Mathematics at Princeton in 1979 with advisor Fred Almgren, and, then arrived at Courant Institute, NYU, in the Fall of 1979 as an NSF Postdoctoral Fellow to work with Louis Nirenberg. We return to this  below.

Bob Kohn remained at Courant for his entire career, advancing from
Assistant Professor (1981) to Silver Professor of Mathematics (2017) and Professor Emeritus (2022),   until he passed away.

During his years (1975-1979) as a Ph.D. student with Fred Almgren (a
geometric measure theorist) and an NSF Graduate Research Fellow at Princeton University, we find Bob Kohn writing his innovative thesis on the theory of elastic deformations, deeply influenced by Fritz John \cite{MR2630218}. This research laid a foundation for his lifelong interest in applied mathematics, partial differential equations problems, and in particular in the mechanics and the mathematics of materials (Richard James discusses this aspect of Bob Kohn's work in more detail below).

After completing his Ph.D. in 1979,  Bob Kohn arrived at the Courant Institute, NYU, the citadel of applied mathematics and differential equations. He now was an NSF Postdoctoral Fellow, with the intention of working with Louis Nirenberg. This is a phase transition for Bob.  As told by Louis, returning from a dim sum lunch, Louis observed that he, Bob, and Luis Caffarelli had been spending time together and hence they should work together. Luis Caffarelli was a professor at Courant at that time, working with Nirenberg.
 The question was, on what? Bob, knowing Scheffer's work \cite{MR454426}, suggested the Navier-Stokes equations, and a collaboration was born, which led to a seminal contribution to the existence and partial regularity results of the weak solutions of the Navier-Stokes equations \cite{MR673830}. This work of Luis Caffarelli, Bob Kohn, and Louis Nirenberg, and in particular, the partial regularity result, which is known as the Caffarelli-Kohn-Nirenberg Theorem, has become a foundational in mathematical fluid mechanics and was recognized by the 2014 AMS Leroy P. Steele Prize for Seminal Contribution to Research, \cite{SteeleP}. 
Understanding the behavior of solutions of the Navier-Stokes Equations remains a major unresolved area of fluid dynamics. It is listed among the still unsolved Clay Millennium Prize Problems. The  Caffarelli-Kohn-Nirenberg work is an important part of the formulation of this celebrated Clay problem \cite{ClayP}! (by C. Fefferman)
As discussed in detail below by Bob's close collaborator and friend Yoshikazu Giga, the famous Caffarelli-Kohn-Nirenberg work motivated their long-term collaboration too, first with the goal to construct a blow-up solution for the Navier-Stokes equations. In their quest (see Giga's fascinating recollection below for full details of the collaboration), Bob and Yoshikazu started with simpler models of semi-linear parabolic equations. This exploration led to another series of groundbreaking results on the asymptotic behavior of  solutions near the blow-up point for such problems \cite{MR784476,MR876989,MR1003437} as well as opened new directions in analysis and differential geometry. Further, we note that the 2026 Breakthrough prize in Mathematics is directly connected to the above research of Bob (with Caffarelli, Nirenberg and Giga). It was awarded to Frank Merle  ``For breakthroughs in nonlinear evolution equations, with regard to their stability, singularity formation, or resolution into solitons'',  \cite{BRP}.
\begin{figure}
\centering{\includegraphics[width=0.75\textwidth]{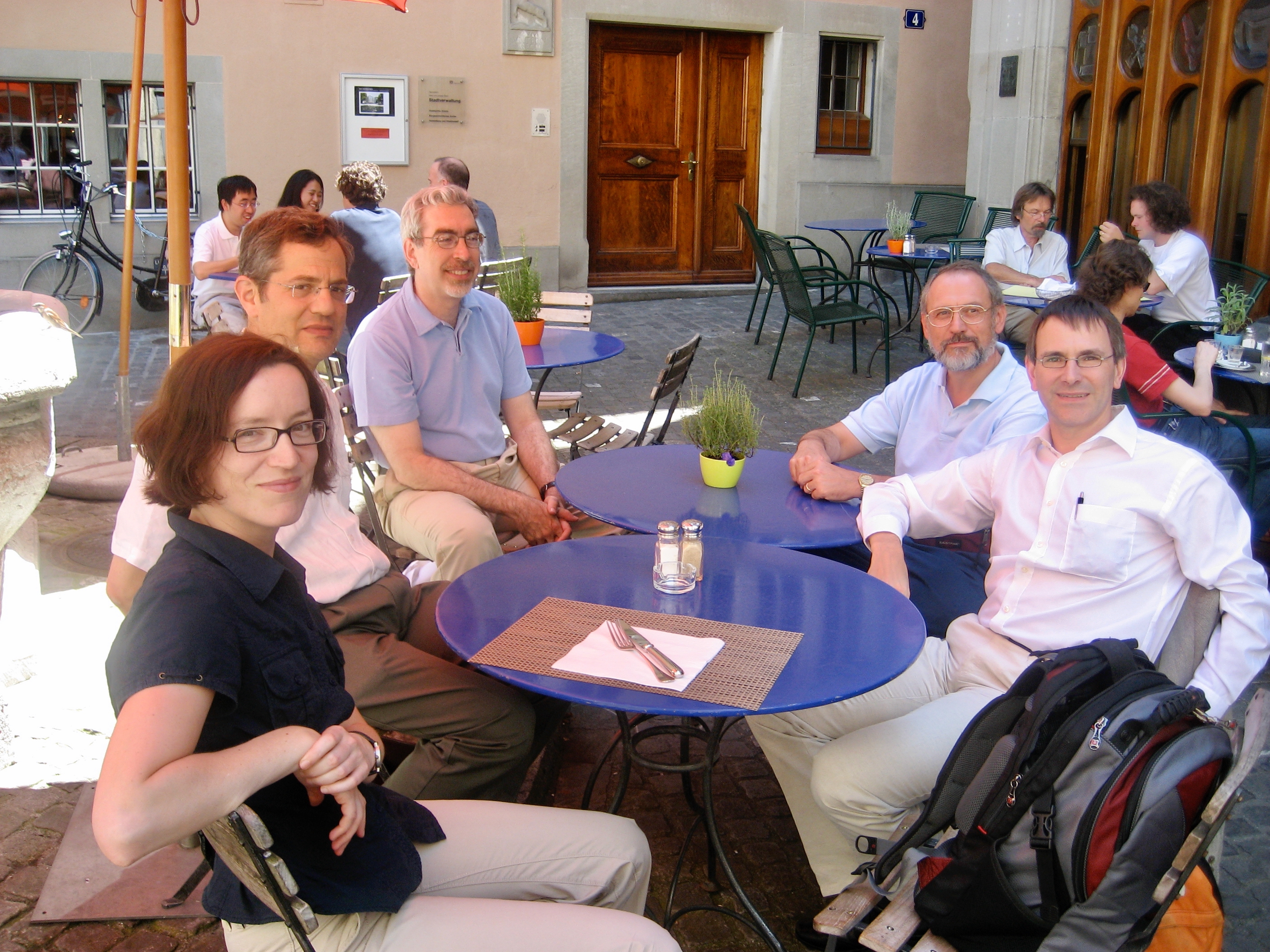}
    \quad
\includegraphics[width=0.75\textwidth]{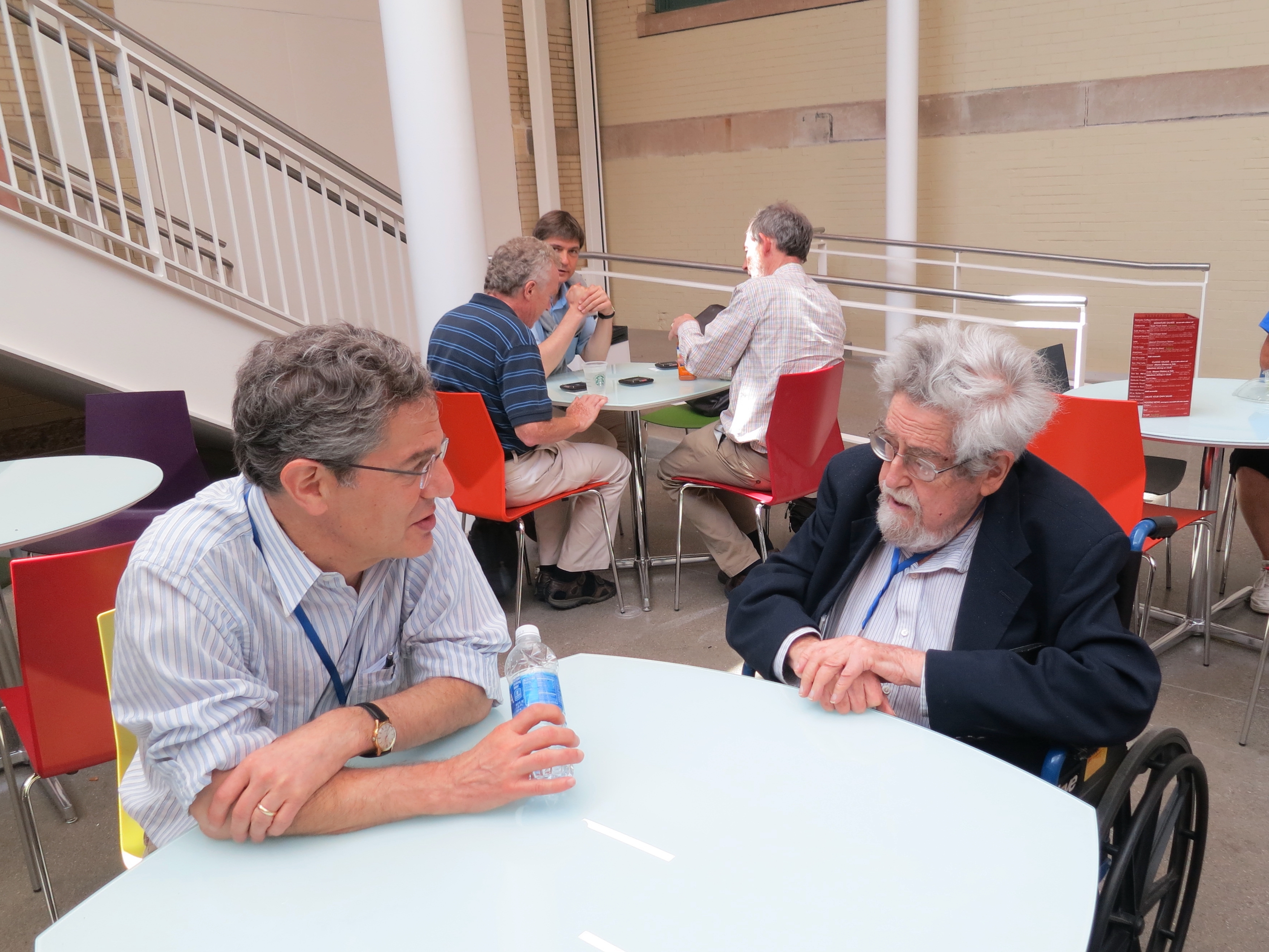}}
\caption{\footnotesize R.V. Kohn at ICIAM meeting in Zurich, 2007
  (top). R.V. Kohn and L. Nirenberg in Pittsburgh, 2016
  (bottom). Photos due to D. Kinderlehrer.
} \label{fig2}. 
\end{figure}  

Bob's mathematical achievements were recognized once again with a Sloan Research Foundation Fellowship (1984-1986). During these years Bob 
started to transition his major research focus back to
 emergent issues in mechanics and mathematics
of materials, the domain where his deep and fundamental achievements would be realized throughout
his entire career. And here he would become the
leading figure in applied mathematics, defining and
promoting major research directions.

This is also the area where Bob would serve as a remarkable mentor, generous and inspiring collaborator to generations of young scientists. Bob Kohn mentored and shared his knowledge and wisdom with many students. He advised 36 Ph.D. students of his own, and mentored so many junior researchers and postdoctoral associates. Bob's keen interest in the field of mechanics and mathematics of materials focused on research areas ranging
from homogenization, optimal design of composites and structures, phase
transformations, physical systems with surface-energy-driven coarsening and their coarsening rates, plasticity, plate theory, and pattern formation to
electromagnetism, micromagnetics, superconductivity, among other
topics as discussed in detail below by Antonio DeSimone, Irene Fonseca, Richard James, Stefan M{\"u}lller, Felix Otto and Michael Vogelius. As a part of
this absolutely fascinating work, Bob established collaborations with Gilbert Strang, Roger Temam,
Michael Vogelius, Michael Weinstein, among many other researchers. Moreover, in this area of research, Bob Kohn 
mentored many postdoctoral fellows, among them are Kaushik Bhattacharya, Stefan
M{\"u}lller and Felix Otto, as well as advised talented Ph.D. students, among them are Lia
Bronsard, Ian Tobasco and Maria Westdickenberg who provided below their memories and stories about Bob Kohn
and their collaboration with Bob. His former postdoctoral mentees and
Ph.D. students, all wholeheartedly emphasize Bob's amazing enthusiasm, encouragement,
generosity, 
warmth and care for young scientists.
\par Ever conscious of his responsibilities, Bob was profoundly
invested in organizational and educational programs for the
community. Perhaps this began as early as 1981 with the preparations
for the 1984-1985 Institute for Mathematics and its Applications (IMA) theme year ``Continuum Physics and Partial
Differential Equations''. Jerry Ericksen and David Kinderlehrer were
the principal coordinators. In those pre-Zoom years, when we had only
post-it notes, Bob traveled to Minneapolis multiple times to
collaborate. He might say, ``let's hear from the homogenization
people'', and a workshop with opposing points of view would magically
be scheduled. Sometimes such workshops/discussions would occur in the
evenings,  with beer for Jerry, wine for David, and hot water with
lemon for Bob. It is noteworthy to mention that Bob served on the
Board of Governors for IMA (2004-2008) and was the Chair of the
Organizing Committee for the Fifth Abel Conference, ``Celebrating the Mathematical Impact of John F. Nash Jr and Louis Nirenberg'' (2015) at IMA.

There are many meetings and conferences where we remember Bob's
role. A highpoint is that he, with David, organized the first SIAM
conference on mathematical aspects of materials in 1994 and then the
second in 1997. They reprised this in 2010 with the first SIAM conference of the newly organized SIAM Activity Group, of which Bob was the chair. Another highlight is Bob's engagement in the scientific program of the ``CMDS: International Conference, Continuum
Models Discrete Systems'' series and his active involvement and participation in ``CMDS--13'' in  Salt Lake City in July 2014 - a highly interdisciplinary conference that brought together scientists working on continuum theories of discrete mechanical and thermodynamical systems in the fields of mathematics, theoretical and applied mechanics, physics, materials science and engineering. The above are just a few illustrative examples as Bob was involved in organization of many other scientific conferences, workshops and meetings all of which benefited significantly from Bob's scientific vision and insights.
\par Bob played a crucial role at Courant Institute, NYU,  his home institution too, where he took part in major activities, including being fully engaged and committed to the success of the Courant/NYU ``Financial Mathematics Master's Program'', as well as being an integral member of the Institute's Steering Committee for many years as discussed by his close colleague and friend Russel Caflisch below. Bob was also a supportive and generous colleague to many as communicated by Leslie Greengard below, as well as recognized by various events that Bob Kohn organized at Courant. For example, most significantly,  Bob organized ``The Memorial in Honor of Louis Nirenberg'', with Jalal Shatah, January 19, 2023, at Courant. Furthermore, Bob served SIAM community for more than 40 years in various roles. Bob's extraordinary scholarly achievements were honored in the greater applied mathematics and
education community by his election to the American Academy of Arts
and Sciences (2017), as Fellows of AMS (2012) and SIAM (2009), via
several plenary lectures, including by Bob's ICM Plenary Lecture
``Energy-driven pattern formation'' at the International Congress of Mathematicians (ICM) in Madrid in 2006 \cite{MR2334197}, by Keith Medal
(Royal Society of Edinburgh, 2007, shared with A. DeSimone, S. M{\"u}lller, and F. Otto for work \cite{MR1854999}) and Ralph
E. Kleinman Prize (SIAM, 1999) among various other distinctions.
\par On a personal level, Bob and his beloved wife, Leslie Anker, met
on a cross-country ski trip in Vermont. They shared a profound passion
for arts, travel, people, nature and sports, and had countless
memorable experiences.  Bob enjoyed cooking and entertaining, and
their home was always open to family, students, friends, and
colleagues.  As Leslie confirmed, in all aspects of Bob's life, he
always aimed to be the best and encouraged others to achieve their
goals. Bob Kohn will be deeply missed by so many friends, colleagues and students who had the great privilege and honor to know him. Our warmest thoughts and heartfelt sympathies go out to Bob's beloved wife Leslie Anker and their family.
\par Below we present stories and memories of some of Bob Kohn's friends, colleagues, postdocs and students.
\section*{Memories}
{\bf Kaushik Bhattacharya (Caltech):}\\  Perception.  Bob had an incredible ability to rapidly grasp the core issue and strip out the unnecessary details when faced with a complex problem.   In the fall of 1990 when I was a graduate student in Minnesota, I found myself on a hike with Bob.  Like he always did when he saw a student or young researcher, he walked up to me and asked me about my work.  Excited, I launched into a long explanation of my thesis.  As I finally ran out of steam, Bob looked at me and said, ``Listen, let me see if I understand'', and proceeded to summarize everything I had said in one concise sentence.  I was stunned: he had not only listened to me, but had grasped the essence and said all that needed to be said.  The next week, at a workshop, I saw him do this to many seasoned researchers.  
 
Educator.  Bob always gave his attention to, and encouraged young researchers.  As I noted above, he  walked up to graduate students and postdoctoral scholars in a mixed group and made them feel like a part of the community.  He seemed to know every student at Courant, and often stopped by the student's lunch tables on the 13th floor of Courant.  And he was indefatigable in visiting student posters at workshops and conferences.  These conversations were not perfunctory, but meaningful conversations that left an impression on the young researcher.
 
Mathematics.  Bob was able to make connections between quite diverse
fields.  In the study of phase transformations, Bob made a remarkable
connection between the Calculus of Variations approach of
Ericksen/Ball-James/Chipot-Kinderlehrer, and the physics theories of
Roytburd-Khachaturyan-Sethna, \cite{ISI:A1988Q342600002,MR906132,MR556558,osti_5821133,ROITBURD1978317,osti_10169268}.  Further, he introduced ideas from the
study of composites (translation method and compensated compactness)
to the study of phase transformation.\\
{\bf Lia Bronsard (McMaster University):}\\
Here is a sampling of the many, many fond memories I have with Bob: I was Bob's 4th doctoral student, and I remember going to his door (which was always open) in the process of choosing a thesis advisor, to discuss potential projects. He gave me two suggestions, which I found so interesting that my choice was made. Both projects led to fundamental new directions in mathematics, the one that I chose leading to the many beautiful results on mean curvature flows and their relation to the Allen-Cahn functional in materials science \cite{MR1075075,MR1101239}. Later, during my post-doc at the IAS, I went to visit Bob and once again he proposed a beautiful new direction on vector-valued Allen-Cahn problems, which led to the first results on the short-time existence for the evolution of networks of curves by their mean curvature \cite{MR1240580}. Both topics have led to a huge literature and beautiful and profound analytical and numerical results in the mathematics connected to material sciences, and formed the foundation of my research career.
This would continue and over the past 25 years, every time I met him, he would listen to my results and directions and would give me suggestions, writing down his ideas in detail and giving me many references related to the topics under discussions. I will always remember his huge cabinets full of articles, and him knowing exactly what to look for, given the discussion we were having. His time seemed endless as he was sharing it: there was no pressure to leave his office. I know I was not alone in receiving Bob's largesse; his openness and generosity were there for generations of students and postdocs.
In the last ten years or so, we would go for lunch when I was visiting him, and this has also allowed me to see his very thoughtful and generous personality. He would also gladly share his time on e-mail: whenever I would ask his advice, for example in recommending plenary speakers or for important issues to address as a SIAM Activity group chair, he would not only give many suggestions, but also detailed explanations behind each. It has been very inspirational to see how open he was in sharing his valuable time with everyone equally. 

My last meeting with him was perhaps the most touching: he asked me if I would be in NY in the Fall 2025 and so I went to meet him during my fall break. He spent all afternoon with me, having lunch, discussing many aspects of our lives as mathematicians, and discussing all my projects. He was again deeply interested, giving me names of articles related to each project, and suggestions for new directions. One of them relates my original doctoral work to new fractional partitioning problems, and I am so happy that I had time to tell him that his question related to heteroclinics for fractional Allen-Cahn was now resolved! 
Great mathematicians that accept to share their insights with all of us are perhaps the greatest of all, and I am very fortunate to have been a student of such a great mathematician.\\
{\bf Russel Caflisch (Courant Institute, NYU):}\\ It is my sad but great privilege to share a few memories of Robert Kohn, a cherished colleague from our time at New York University's Courant Institute of Mathematical Sciences. I've known Bob since the 1980's when we were both young faculty members at Courant, but we had never worked closely together until I returned to Courant to serve as Director in 2017. Bob was an integral member of the Institute's Steering Committee, and I came to rely on his character and sound judgment throughout my term.

Bob brought a steady, clear-eyed perspective to every discussion. He had an exceptional ability to understand complex administrative issues and balance differing viewpoints. This was a natural extension of Bob's mathematical skill, where he was remarkable for both the depth of his scientific insight and his ability to connect ideas and people across disciplines. He envisioned and enabled the application of mathematics across a wide variety of fields.

Bob's wisdom, integrity, and quiet leadership left a lasting impression on myself and on the entire Courant community. He will be deeply missed by generations of students, colleagues, and friends.\\
{\bf Antonio DeSimone (SISSA):}\\
My first strong memory of Prof. Robert V. Kohn is from 1992. It is from a summer conference in Edinburgh, 
where I met him in his role of post-doctoral mentor of a friend and former fellow-graduate student K. Bhattacharya. Seeing him in action  convinced me immediately that I really wanted to have an experience as a postdoc, even though I already had an academic job.

I then met Bob, in his capacity of informal mentor of whoever was smart enough to seek his advice. Like many others, I experienced several ``RVK-moments''. I would tell him about my current research problem, he would listen silently for a while to my results and to my questions and then politely reply: ``Let me see if I understand''. His version of my problem had suddenly acquired structure, clarity  and relevance way beyond what I had been able to see. My own experience as a mentor owes much to these RVK-moments. The times I manage to do something remotely similar with my students and post-docs are the ones that I remember as the most successful of my career.

I spent quite some time reading, rereading, mastering, and admiring his 1991 paper with the explicit calculation of the relaxation of a double well energy \cite{MR1122017}. It taught me the beauty and the elegance of an explicit result in a very concrete problem of materials science. Very concrete, but relevant to the point of being paradigmatic. Throughout my own ten-year-long research adventure on the relaxation of non-convex energy functionals (the relaxation of the micro-magnetic energy of hard ferromagnetic materials in the large body limit and the quasi-convexification of the free-energy of nematic liquid crystal elastomers, the latter obtained jointly with G. Dolzmann \cite{MR1894590}), his 1991 paper has been a constant source of inspiration.

In the following decade I had the privilege of collaborating with him directly, together with S. M{\"u}lller and F. Otto, on a broad research program on thin ferromagnetic films which has lead to several ``DKMO'' papers, e.g. \cite{MR1916988}.
Among the many things I learned in this period is that even partial results can provide tremendous insight and even be lifted to the status of general ``methods''. This was the case for the approach by matching upper and lower bounds for the study of the scaling behavior of energy functionals with respect to asymptotic regimes of geometric and material parameters, that he pioneered.
In the case of ``DKMO'', this method was applied to the energy of thin ferromagnetic films and its scaling behavior with respect to small parameters, such as (for example) the thickness of the film or the magnitude of the anisotropy constants, but a similar approach has been proved to be very effective in a wide range of problems by an entire generation of researchers.

We, all of Bob's friends, will sorely miss him. But actually, from a personal and possibly selfish perspective, the first things that come to my mind when thinking of Bob are joyful rather than somber: the vivid memory of the many lessons learned from him, the feeling that they will keep accompanying me for the rest of my life, and how lucky I have been to meet him.\\
{\bf Irene Fonseca (Carnegie Mellon University):}\\
What I remember most about Bob was the generosity with which he shared his thoughts and possible open problems with others, especially with young(er) researchers. In the early 1990s, shortly after Luc Tartar and I wrote a paper in 1989 on a singular perturbation problem in the gradient theory of phase transitions \cite{MR985992}, Bob invited me to visit him at the Courant Institute. Bob and Peter Sternberg were also interested in these models, having written a paper in 1989 on local minimizers for these energies and their asymptotic behavior \cite{MR985990}.
We spent a couple of days discussing possible directions to undertake in this context, and when I was leaving, Bob told me that he was sorry we had not had time to begin a collaboration, but that he hoped some of these ideas would be pursued. And they were, in later work with some of my trainees. This was typical of Bob: he was generous in sharing his mathematical creativity, intuition, and insights. He was, and will always be, a source of infinite inspiration.\\
{\bf Yoshikazu Giga (University of Tokyo):}\\
Bob Kohn was an exceptional mathematician who was really a great interdisciplinary scientist. He was strong at organizing substantial discussions with various people, not only mathematicians but also people in other disciplines. As a result, he started many new fields in mathematics related to applied fields like materials science. His mathematics was so fundamental that his theory even has a strong influence on the original applied fields. In mathematics, he often clarified the hidden relation between two apparently different fields. For example, Bob Kohn (with S. Serfaty) found a deep relation between some discrete game (in combinatorics) and motion by its mean curvature (in materials science and differential geometry), \cite{MR2200259}. He was an ideal mathematician representing the spirit of ``Courant Institute of Mathematical Sciences (CIMS)'' and ``Communications on Pure and Applied Mathematics.''
Bob Kohn and I met first in Tokyo in July 1982. He presented his very famous result (obtained jointly with L. Caffarelli and L. Nirenberg \cite{MR673830}) on the partial regularity of suitable weak solutions of the Navier-Stokes equations both at the US-Japan seminar and at Nagoya University, where I was working at that time. Since I was also working on the Navier-Stokes equations, we started a discussion. From November 1982, I stayed at CIMS for one and a half years. Bob and I were trying to figure out if we could construct a blow-up solution (a solution whose value diverges in finite time) for the Navier-Stokes equations (later known as Clay's millennium problem and it is still an open problem). Bob suggested constructing a self-similar (finite energy) blow-up solution following the idea of J. Leray in his famous 1934 paper, \cite{MR1555394}. (A negative answer was given in the 1990s.) Bob also expected that such a solution would describe a typical behavior of blow-up. Of course, the Navier-Stokes equations are difficult. We began to investigate a similar problem for parabolic problems for which we know that there exists a blow-up solution. A simplest example is the Fujita-type semilinear heat equation, a heat equation with power nonlinearity causing blow-up. We chose this equation to start our project. At that time, the focus was on whether a blow-up would occur for a given initial value, and there was little research on the behavior of the blow-up. This is why we began studying blow-up behavior.\\
The first topic we focused on was the existence of non-trivial (not uniform in spatial direction) self-similar blow-up solutions. At first, we tried hard to construct such a solution. To my surprise, Bob tried to look it up on his pocket calculator to study the one-dimensional case. This is a good example of how Bob was not concerned with the means to find the essence of the problem. We got the impression that a non-trivial self-similar solution does not exist in lower dimensions. Based on this observation, we found a Pohozaev-type identity which proves the non-existence of a non-trivial self-similar blow-up solution for subcritical Sobolev exponents. This was an early stage of our joint work, which resulted in three joint papers on the asymptotic behavior near blow-up point \cite{MR784476,MR876989,MR1003437}. Our study produced a lot of new problems, and we sometimes came back to these topics and methods. Fortunately, our naive results and methods led to important progress in various other problems, not only in mathematical analysis but also in other fields like differential geometry. For example, G. Huisken classified the asymptotic behavior of a shrinking mean curvature flow \cite{MR1030675}.

Bob Kohn was particularly fond of intuitive considerations that do not directly lead to rigorous proofs. This is very important for mathematicians to discuss with people in other disciplines. Bob Kohn left a great legacy that strengthens ties between mathematics and all other disciplines.\\
{\bf Leslie Greengard (Courant Institute, NYU):}\\
Bob was an inspiring colleague and an extraordinary mentor. From the time I first arrive at Courant as a postdoc, through our decades together on the faculty, he was one of the first people I would seek out for advice on analysis, on selecting research topics, on exploring new ideas, and on navigating  life as a teacher and advisor. I miss his warmth and wisdom.\\
{\bf Richard James (University of Minnesota):}\\
Kohn's interest in mechanics -- elasticity theory in particular -- began with his thesis, ``New estimates
for deformations in terms on strains''.  It was an interest that followed him through his life, his last paper 
 being ``The effective energy of a lattice metamaterial'' \cite{MR4989947}.  
Within mechanics, his selection of major topics included the regularity of solutions of the Navier-Stokes equations, 
homogenization, optimal design of composites and structures, phase transformations, plasticity, plate
theory, and pattern formation.   He also pursued ongoing research on electromagnetism, micromagnetics, 
superconductivity and finance.
\medskip

Kohn did his graduate work at Princeton University under the direction of Frederick J. Almgren, Jr., whom he
thanks for ``unfailing optimism and confidence''.  His thesis (141 double space pages) is notable for its
clarity.  Kohn explains that ``...our conventions will be those of H. Federer,
\underline{Geometric Measure Theory}'', \cite{MR257325}.   But the thesis also contains a concise chapter on 
``The model of hyperelasticity'' from which one can pretty much gain a working knowledge of elasticity theory 
by reading 7 pages.  Perhaps Kohn's is the only PhD thesis written with the conventions of 
Federer that gives measured material constants of an accepted model of rubber.
\medskip

In the early 1980s at the Courant Institute, homogenization and the optimal design of composites were topics
of intense research, pursued especially by Graeme Milton, Marco Avellaneda and Kohn.   Among many others, 
Fran\c{c}ois Murat and Luc Tartar, Konstantin Lurie and Andrej Cherkaev,
and Gilbert Strang (with Kohn) rapidly developed the subject.  The methods with limited application were
quickly discarded, and those with broad application became linked more closely to  fundamental concepts 
(and hard open problems) in the calculus of variations, and, more broadly, in analysis.  His paper with Gilbert
Strang \cite{MR707959} on optimal design laid out the fundamental obstacle in 1983:
\medskip

``Unfortunately, quasiconvexifications are hard to compute. So far, all examples have involved ordinary 
convexification in an essential way. 
An underlying difficulty is the lack of an algebraic condition for quasiconvexity. ($u:\R^2 \longrightarrow \R^n$ , 
$G(\nabla u)$ is polyconvex if it is a convex function of $\nabla u$ and its $2 \times 2$ minors.) 
Hadamard gave a necessary condition, namely rank-one convexity (also called ellipticity, or the 
Legendre-Hadamard condition). But polyconvexity is not necessary; the sufficiency of rank-one convexity 
is open; and the condition $QG = G$ is neither algebraic 
nor easy to work with.''

\medskip

\noindent (The sufficiency of rank-1 convexity failed as well, from a 1992 example of Vladimir Sverak \cite{MR1149994},
and $QG = G$ is not algebraic \cite{MR1668552}.) 
\medskip
 
The early 1980s were a period of intense activity for Kohn.  Not only was he a  leading figure in homogenization
and optimal design, but, on a  different subject involving completely different mathematical methods,
he proved a ground breaking theorem on partial regularity of the Navier-Stokes equations 
with Luis Caffarelli and Louis Nirenberg (1982, 1984). 

Following his main body of work on optimal design, Kohn set new directions of research.  
He realized that the  connection with quasiconvexity linked homogenization and optimal design with broad 
areas of materials science.  The general paradigm is a variational principle with a nonconvex integrand,
$\varphi (\nabla \bfy, \bfy, \bfx)$, in the vector-valued case.  Such energies are ubiquitous: from structural phase 
transformations and micromagnetics (and its analog for ferroelectrics) to the Ginzburg-Landau theory of 
superconductivity.  With Stefan M{\"u}ller he developed a theory for the branching of twins at an interface 
between phases \cite{MR1272383}.  The theory enjoyed both intense development by mathematicians as well as diverse
new applications in materials science.\\
{\bf Stefan M{\"u}lller (Hausdorff Center for Mathematics, Bonn University):}\\
What strikes me most about Bob,  aside from his being an outstanding mathematician per se, is his enormous talent to read the literature in other fields in depth, to listen to and engage with people in all areas, as well as to identify and isolate problems of fundamental interest in applications where mathematics can really make a difference in \emph{ understanding}. Almost every time I saw him he had identified a new problem in the science literature and put it in a form which led to fascinating mathematical questions. Of course, the other striking feature, which almost everybody at Bob's memorial service alluded to was Bob's scientific and personal generosity. Now journals are awash with articles that reproduce experimental curves without regard to the actual phenomena. This was not Bob's way. Bob's understanding included the fact that a theory should reproduce observations (or predict new ones) for the right reason. Bob's approach to applied mathematics had a strong influence on me and in fact I think on a whole generation of applied mathematicians. I cannot think of anyone who so deeply tried to understand other people's work and way of thinking (in science or in mathematics) while keeping his own very original viewpoint.

An example is Bob's work with Kaushik Bhattacharya about the recoverable strain in polycrystals  \cite{BHATTACHARYA1996529,MR1478776}. The known shape memory methods based on reversible single crystal kinematics did not work. They employed with great success an original extension of the duality theory in the form of Voigt and Reuss bounds for mixtures that accommodates the polycrystalline structure. Other applications include inverse problems, e.g.,  \cite{MR739921}  (with Vogelius, predating Bob's work on materials), cloaking, e.g.,  \cite{MR2384775}, coarsening rates \cite{MR1924360} (with Otto), branching of microstructures, micromagnetics, e.g.,  \cite{MR1272383,MR1854999} (with DeSimone, Otto, me and others),  folds/blistering/wrinkling/draping, e.g. \cite{MR985990,MR1752602,MR1854999,MR3179665,MR3628880,MR4215193}  (with Jin, Sternberg, Bella, ..., later continued by his student Tobasco), crystal growth, and of course a whole body of work on composites and optimal design.

`` ... we seek understanding as well as answers'' is a quote from the introduction to a survey article on micromagnetics, jointly with Antonio DeSimone, Felix Otto and me. In this case the insight that micromagnetics is an area ripe from mathematical investigation came from Antonio. Bob and I joined later and finally Felix took the analysis to a new level. The article is a report on joint work. As you can imagine, the introduction was written by Bob and in my view summarizes very nicely his thinking.
This survey appeared in a three-volume work on ``The science of hysteresis'' (Eds. Giorgio Bertotti and Isaac D. Mayergoyz).

I first met Bob in April 1986 in Oberwolfach at a meeting on Elasticity Theory organized by John Ball and Willi J{\"a}ger, but we really started to work together when I was a Zeev Nehari Asst. Prof. at CMU and Bob invited me to the Courant Institute in March 1990. This was the beginning of a long and very inspiring collaboration starting with branching in a simplified model of an austenite/finely twinned martensite boundary (
let me mention that as a young German mathematician I was very moved by the very warm welcome I received both at CMU and at Courant). 

Regarding historic developments, I sometimes asked myself what made Bob move from great successes in more ``pure'' analysis (Cafarelli-Kohn-Nirenberg \cite{MR673830}, Giga-Kohn \cite{MR784476,MR876989,MR1003437}) to his personal approach to ``applied'' mathematics. Perhaps, as Michael Weinstein intimated at the Memorial Service, it had always been in Bob's mind to understand the world through mathematics, as Bob himself wrote , quoted above, ``... we seek understanding as well as answers.''\\
{\bf Felix Otto (Max Planck Institute - Leipzig):}\\
To my generation, Bob was the good soul of the Courant Institute, in the tradition of Peter Lax. In particular, he mastered the art of bringing people from different scientific backgrounds together, always with a fruitful and selfless plan in mind. As a necessary first step towards this, he loved to introduce people, which he did in his charming and casual way. Occasionally, he introduced people who knew each other already pretty well - which of course he would instantly realize: ``I guess you don't need introduction'' with a disarming smile. 

I was asked to comment on our joint 2002 paper ``Upper bounds on coarsening rates'' \cite{MR1924360}. While my anecdotal memory is pretty bad, I clearly recall that it was an easy - and short - paper, in the sense that it was conceived and mathematically pinned down within a couple of weeks, while I was Bob's guest at Courant. From the applied mathematics endeavor on liquid crystal, which Bob co-initiated, we were familiar with spinodal decomposition and the coarsening exponents describing its later stages, as documented in a rich experimental literature. We were also familiar with some theory, in particular in the context of Ostwald ripening (via Barbara Niethammer, a former Courant postdoc mentored by Bob) and grain growth (via David Kinderlehrer). In terms of partial differential equations, it was the high time of gradient flows and energy-energy dissipation relations, and both of us had leveraged these structures, Bob in his work with Yoshikazu Giga \cite{MR784476}. The key concept was a notion of an average length scale of the phase distribution that is strong enough to control the interfacial energy from below, but weak enough to be dynamically controlled by energy dissipation. It turned out that a transportation distance between the two phases was flexible enough to deal with both shallow and deep quench regimes. It was a bit of a mental hurdle to accept the fact that the outcome is a time-averaged instead of a pointwise upper bound. The only slightly technical ingredient of the proof turned out to be an interpolation estimate Bob and I had come up with in a 1999 paper with Rustum Choksi on domain branching in uniaxial ferromagnets \cite{MR1669433}. Bob had realized that the mathematical nature of the driving forces behind domain branching in ferromagnets was pretty much the same as of those that made twins split in Martensites, a mechanism he pinned down in his seminal 1994 work with Stefan M{\"u}lller \cite{MR1272383}. Loosely speaking, time in the coarsening problems plays the role of the longitudinal axis in the materials science problems, and the same interpolation estimate encapsulates the competing forces. Incidentally, the introduction of a transportation distance proved valuable when dealing with flux-tube branching in type-I superconductors in a work with Sergio Conti, Sylvia Serfaty and myself \cite{MR3542007}. The range of physics phenomena captured and the unideological use of diverse mathematical tools in this sequence of works is typical for Bob's oeuvre.\\
{\bf Ian Tobasco (University of Michigan):}\\
Bob Kohn was my PhD advisor from 2011-2016, and I count myself lucky to have experienced his mentorship firsthand. I first met Bob as an undergraduate who had applied to the PhD program at Courant. Bob emailed me to encourage me to accept their offer and volunteered to speak by phone. When I finally arrived in New York, it was on the same night as hurricane Irene. Bob and his wife Leslie opened up their home to me to have company and shelter from the storm. Afterwards, they sent me home with a grab-bag of goodies, so that I wouldn't have to fuss over groceries, at least for a little while.
Over the years, I learned just how much Bob loved being a mentor. During the PhD, Bob met with me at least once per week---sometimes two or three times---either telling me to come on in to his office when I stopped by, or taking out his pocket book to pencil in another time. Meetings were both therapeutic and lively. One time, I explained a mess of calculations that seemed to go nowhere involving the Jin-Kohn entropy method \cite{MR1752602} and a problem in elasticity. Bob smiled and said, ``Well, you learned something.'' Another time, I was supposed to have prepared to present a paper of his on the shape of leaves and flowers, but came with notes about crumpled cylinders instead. Bob simply pivoted and asked me to explain what I was thinking about. Bob brought a zen-like atmosphere to doing mathematics and taught me that if one thing didn't work out, then another could. Bob was the model of an even-keeled and deep researcher who knew how to, as he liked to say, ``play the long game.''
Bob was extremely generous with his ideas, to the point that he was unwilling to take credit as a co-author unless he was convinced that he had done something new. Just ``being a sounding board'' did not count. Once, I was trying to wrap up a paper that would eventually be part of my thesis, but didn't know what to say in the introduction. I emailed Bob asking for advice, and the next day found the first half of the introduction in my inbox, along with a point-by-point explanation of the value of each paragraph. When I asked him to be a co-author, he refused, saying that I had done all the hard work.
Bob kick-started my academic career and helped me find my way. He demonstrated to me time and time again how to ask good mathematical questions, how to reach across disciplines to build bridges, and how to be a caring mentor. When I find myself stuck in the weeds, I often think, ``What would Bob do?'' I know I'm not the only one.\\
{\bf Michael Vogelius (Rutgers University):}\\
Bob and I met for the first time in early 1981. We were then both Postdocs at NYU, and George Papanicolaou suggested we might have common interests. As it turned out, he was very right. A little after our first meeting Bob showed up with a paper he had just received from David Jerison. It was a small 1980 note written by Alberto Calderon, and the rest is history: that was the beginning of our first collaboration. Fast forward to 2007 - Bob and Mike Weinstein had just been out here at Rutgers, and a couple of days later Bob calls me and asks, if I have seen the Science Section of the NY Times, there was an article about ``cloaking''- and this became the topic of our last significant collaboration, now also involving Mike, who incidentally had been a student with George, when we were Postdocs at NYU, e.g., \cite{MR2384775}. As both of these incidents show: in addition to being a formidable mathematician, Bob was always very sharing and enthusiastic about including others in his discoveries. 
Over the years Bob, Mike and I became close personal friends and not just collaborators. Indeed, after a very nice meeting at Luminy, where our wives, Leslie, Suzanne and Kirsten also got to know each other well, we started a tradition of alternating lunches in NYC and rural NJ. We will try to continue this tradition, but it won't be the same without Bob. I will personally miss him a lot, and so will the entire community.\\
{\bf Maria Westdickenberg (RWTH Aachen University):}\\
Bob Kohn had a gift for seeing straight to the nature of a problem. He had an amazing gift for connections and a remarkable generosity for sharing those connections with colleagues, students, and all the many guests who visited his office on a trip to or through New York. As deep and as rich as many of his papers were, shedding new light on singularities and coarsening and energy-driven pattern formation and so many intriguing phenomena from materials science, his penchant for connecting people and ideas is an extra distinguishing feature. 
Probably everyone who ever sat in Bob's office remembers the times he quietly stood, walked over to his filing cabinet, and cooly extracted a paper that was eerily spot on for something one was trying to understand. He gave many impressive talks and hosted many minisymposia, but he said that the most important part of conferences were the conversations that took place in between the talks. He said that the most important talks to attend at one's home institution were the ones by locals, because you should really know what the people around you were working on. He would pull you aside at a workshop and introduce you to someone you really should meet because of a scientific connection he had not only noticed but remembered. He did this not with a narrow aim but because he knew it was good (essential) for the community.
He was also classy. And timeless. He could say things like ``gee'' and ``gosh'' and get away with it because of the insight of his mathematical observations. But he had standards. One should write ``different from,'' not ``different than.'' Language matters; details matter.
Bob was not pretentious. He was always willing to start with something simple in order to get the idea across. Of course he was also willing to continue to arbitrary depths and complexity-- but in the beginning, he wanted to convey the essence and make abstract ideas from analysis personal and vivid. He was able to communicate effectively with a broad variety of scientists - physicists and materials scientists and analysts and numericists and others, young and old - because he knew the language of different communities and cared enough to make the effort.
He spoke as if the problem itself had desires... a structure, and a wish to be understood. There are people who think that if you want mathematics to be beautiful, you should study ``pure mathematics,'' and that if you want to work on applied problems, your task is to find the right answer in whatever messy way that you can. Bob's work in applied analysis proves that it is not only possible to find beauty and structure in applied problems, it is even useful to follow the thread of the problem to its natural core.

\bibliographystyle{siam}
\bibliography{ExampleRefs}

\end{document}